\documentclass{pasj01}
\usepackage[section]{placeins}
\makeatletter
\AtBeginDocument{%
  \expandafter\renewcommand\expandafter\subsection\expandafter{%
    \expandafter\@fb@secFB\subsection
  }%
}

\usepackage[square,sort,comma,numbers]{natbib}

\begin{document}

\title{ {Learning Mathematics through incorrect Problems}}
\author{Hugo Caerols-Palma$^1$, Katia Vogt-Geisse${^1}^*$}%
\altaffiltext{}{$^1$Facultad de Ingenier\'{\i}a y Ciencias\unskip, Universidad Adolfo Ib\'{a}{\~{n}}ez \unskip, Diagonal Las Torres 2640, Pe{\~{n}}alol{\'e}n \unskip, Santiago \unskip, 7941169\unskip, Chile}
\email{katia.vogt@uai.cl}

\KeyWords{Education in Engineering-- Higher Education-- Incorrect Problems--  Mathematics Education-- Teaching Techniques}

\maketitle

\begin{abstract}
In this article we describe special type of mathematical problems that may help  develop teaching methods that motivate students to explore patterns, formulate conjectures and find solutions without only memorizing formulas and procedures. These are problems that either their solutions do not make sense in a real life context and/or provide a contradiction during the solution process. In this article we call these problems \textit{incorrect} problems. We show several examples that can be applied in undergraduate mathematics courses and provide possible ways these examples can be used to motivate critical mathematical thinking. We also discuss the results of exposing a group of 168 undergraduate engineering students to an \textit{incorrect} problem in a Differential Equations course. This experience provided us with important insight on how well prepared  our students are to ``out of the box" thinking and on the importance of previous mathematical skills in order to master further mathematical analysis to solve such a problem.   
\end{abstract}

\section{Introduction}\label{intro}

Solving exercises and problems is fundamental for learning Mathematics and are the base for  teaching this subject. According to Halmos et. al,  problem- solving is the hearth of Mathematics and of a Mathematician's work \cite{HAL}.
Teaching mathematics oriented towards problem-solving and developing strategies that help students improve their problem-solving abilities has been of great importance since the 1970s \cite{POL}.
 
 Schoenfeld et al. in \cite{SCH} explain the importance that problem-solving methods have had to teach mathematics in the last decades of the 20th century. The article comments on  the main focus of learning based on problem-solving while developing mathematical thinking.
The authors suggest to make changes in the curriculum and in the teaching methods in order to develop the students' mathematical thinking abilities. They  propose that teachers assign the following tasks to students:   

\begin{enumerate}
\item Exploring patterns, instead of memorizing formulas;
\item Formulating conjectures, not only doing exercises;
\item Finding solutions, instead of memorizing exercises or procedures.
\end{enumerate}

 The following is an example of a problem that illustrates a way to explain how to achieve the previously listed tasks: 
 
 \textit{Given two initial values $x_0=a, x_1=b,$  determine if the mathematical sequence } 
\begin{equation}\label{rec}
x_{n+2}=\frac{x_{n+1}+x_n}{2}, \; n\geq 0
\end{equation} 
 \textit{converges. }

To explore patterns, a student could study the convergence of the sequence by creating a table starting with the initial values and obtaining $x_2, x_3, x_4,...$ recursively, and this way exploring to which value the sequence converges. Students could then formulate a conjecture of the value of convergence that depends on the initial values chosen. Then, the conjecture should be proven. How to show the conjecture depends on the knowledge of the student and the mathematical tools each student manages, these problems can be approached differently and may lead to distinct proofs. For instance, a student could  prove the convergence using an inductive reasoning, whereas a student who knows how to deduce a direct formula for this type of recurrences can find a more direct and different way to prove it. 

In this article we describe a special type of problems that may help to develop a teaching method that makes students explore patterns, formulate conjectures and find solutions without only memorizing formulas and procedures. These are problems  that when solving them generally their solutions do not make sense in a real life context or provide a contradiction. In this article we will call those problems \textit{incorrect} problems and will provide and discuss a few examples.

 We will start presenting a simple \textit{incorrect} problem that could be used in a basic Algebra course during the first years of undergraduate math studies when studying the Field of Real Numbers or even at a high school level mathematics course. Consider the following chain of reasoning:
 
\textit{Suppose that $a=b$ therefore $a^2=ab$ then $a^2-b^2=ab-b^2$ from the previous $(a+b)(a-b)=b(a-b)$. Dividing by $(a-b)$ on both sides of the equation we obtain $a+b=b$, but since $a=b$ we have $2a=a$ dividing by $a$, gives us  $2=1$. }

The first time students are confronted with this problem they usually feel puzzled when realizing that they have reached a contradiction: $2=1$ is definitely not right. By proposing these types of exercises, professors expect that their students become more careful in each step when using deductions to solve them. Some, who at some point have faced this problem and found the reason for the contradiction seldom forget it: there are dangers when dividing by zero. Exercises such as the one just shown may be used to motivate and explain the importance of a deeper understanding of the implied axioms underlying a theory, in this case  the axioms of the Field of Real Numbers, as Lakatos et al. discuss  in \cite{LAK}.  


The objectives of this work are: first, to detect learning opportunities obtained from analyzing the solutions of \textit{incorrect} problems; second, to explore how a group of students-- that are used to be confronted with direct ``correct" problems-- perform when being exposed to an \textit{incorrect} problem; and third, to discuss the initial perception of the professors-- which was that the use of incorrect problems will have positive effects on learning outcomes-- is confirmed in the case of the specific example and group of students analyzed, aiming to give   important insight which may help when dealing with similar exercises in different undergraduate mathematics courses.
 
In Sec. \ref{ExInc} we  present and discuss examples of \textit{incorrect} problems. Particularly, in Subsec. \ref{LProb} we give an example that has appeared naturally while teaching the course Multivariable Calculus and in Subsec. \ref{ProbMeth} and \ref{Cooling} we provide two problems that are relevant in a Differential Equations course for undergraduate engineering students. In Sec. \ref{Students} we show the results of assigning an \textit{incorrect} problem that uses Newton's Law of cooling to 168 undergraduate engineering students in 7 sections of the Differential Equations course at Universidad Adolfo Ib\'a\~nez in Santiago, Chile. Finally, in Sec. \ref{Concl} we will present a discussion regarding this experience, which we believe will be of interest of those who are always searching for a better way of teaching mathematics and have used or are considering using these type of strategies in their courses.

\section{Examples of \textit{incorrect} problems}\label{ExInc}

\subsection{Motivating the Implicit Function Theorem from a limit problem in two variables}\label{LProb}

Limit problems in Multivariable  Calculus  are naturally useful to analyze and establish conjectures. One interesting example is studying the following limit 
\begin{equation}
\lim_{(x,y)\to (0,0)}\frac{xy}{x+y}.
\end{equation}
These problems are often approached choosing trajectories that pass through the origin $(0,0)$ and using the theorem that states that if the limit of a function along one trajectory differs from the limit on another trajectory, then the limit does not exist \cite{Stewart}.  Hence, probably the first approach, which is standard, used by students will be to consider straight lines through the origin as trajectories, resulting in the limit to be zero, and then study the limit along other  trajectory such as  parabolas or square root functions that pass through the origin, also resulting in the limit to be zero. Then,  sensing that the limit exists and is zero-- since different trajectories result in the same limit--  they may try polar coordinates $x=r\cos(\alpha)$, $y=r \sin(\alpha)$, $0\leq\alpha<2\pi$, $r\geq 0$  in the attempt to prove it, obtaining  
\begin{equation}
\lim_{r\to 0} \frac{r\sin(\alpha)\cos(\alpha)}{\sin(\alpha)+\cos(\alpha)}.
\end{equation}
 In what follows we have a first learning opportunity by identifying a mistake which students often make. The mistake lies in concluding that the limit is zero since $r$ approaches zero without considering that the function 
 \begin{equation}\label{galpha}
 g(\alpha)=\frac{\sin(\alpha)\cos(\alpha)}{\sin(\alpha)+\cos(\alpha)}
 \end{equation}
  has to be bounded for $0\leq \alpha<2\pi$ in order to argue like that. 

\begin{figure}
\begin{center}
\includegraphics[width=0.46\textwidth]{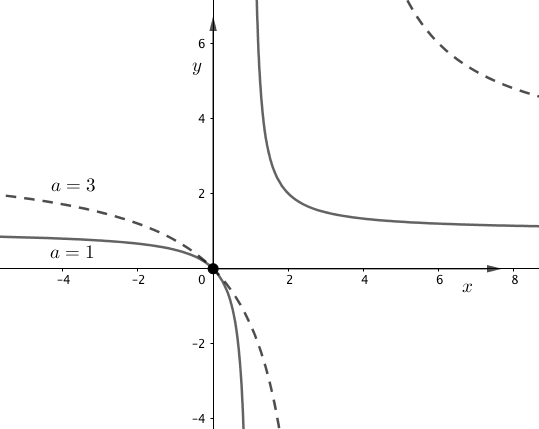}
\caption{Trajectories $\displaystyle y=\frac{ax}{x-a}$ for of $a=1$ and $a=3$.}
\label{figure:1}  
\end{center}
\end{figure}

To make students realize their mistake, one has to prove that the limit does not exist by for example presenting the idea of level curves to students. This leads to an interesting technique to solve the problem. If we consider 
\begin{equation}
\frac{xy}{x+y}=a,
\end{equation}
 solving for $y$  we obtain the curve $\displaystyle y=\frac{ax}{x-a}$, which represents a trajectory that for every $a$ value passes through  the origin and  the limit on that trajectory is equal to $a$. Choosing $a $ different from zero we found a trajectory that gives us a limit different from the zero limit that we obtained along a straight line, this shows then that the limit does not exist. Fig. \ref{figure:1} depicts different trajectories for different values of $a.$ Once proven that the limit does not exists, their mistake can be discussed once again to emphasize the importance to check boundedness of the function $g(\alpha)$ given by equation (\ref{galpha}).

\FloatBarrier

\begin{figure}
\begin{center}
\includegraphics[width=0.46\textwidth]{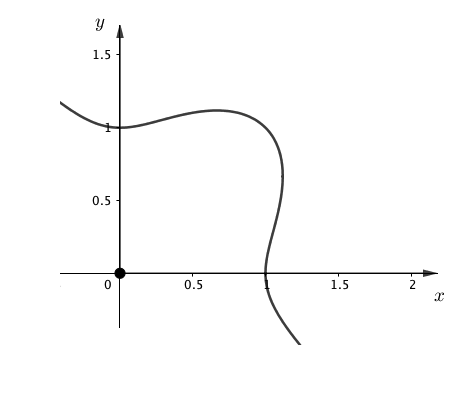}
\caption{Graph given by the equation $x^3+y^3-x^2-y^2=0.$ }
\label{figure:no curva}  
\end{center}
\end{figure}

Afterwards, the students may be  asked to analyze the limit 
\begin{equation}
 \lim_{(x,y)\to (0,0)}\frac{x^3+y^3}{x^2+y^2}.
\end{equation}
 using level curves. 
Equaling to a constant $a$, we obtain level curves of the form 
\begin{equation}
\displaystyle\frac{x^3+y^3}{x^2+y^2}=a,
\end{equation}
 and hence, for $a=1$  we obtain the equation 
\begin{equation}\label{eqn:noCurva}
 x^3+y^3-x^2-y^2=0.
\end{equation}
This is an equation that supposedly represents a curve passing through $(0,0)$ and hence we could argue that the limit along that trajectory is equal to $a=1$, and one can also easily see that the limit is zero along the straight line $y=x$, and hence the argument would conclude that the limit does not exist.  But, observe that by utilizing polar coordinates, it can be easily proven that the given limit is zero. So the question that arises naturally is: What went wrong with the level curves approach? After giving the question some thought, students could observe that the reason the approach went wrong may be  that the relation in Eqn. (\ref{eqn:noCurva}) can not define a function $y=f(x)$ nor $x=f(y)$ near the origin and passing through  it, hence could not represent a trajectory. Indeed, see in Fig. \ref{figure:no curva} the graph of the curve given by Eqn. (\ref{eqn:noCurva}). The origin satisfies the equation but no other point in its neighborhood does. 

From this \textit{incorrect} problem the following question arises naturally: when does a relation between to variables represent one as a function of the other in the neighborhood of a point that verifies it? Thus, we have reached through incorrect limit solutions with two variables an understanding of the real importance of  the Implicit Function Theorem (see the details of the Theorem in \cite{Rudin}).

\subsection{A problem that questions numerical methods}\label{ProbMeth}

We present in this subsection an example that illustrates interesting questions that arise from an initial value problem, for which the initial condition escapes the interval of definition of the solution. In particular we are going to study how a numerical solution method may fail to give expected answer. 
Consider the following \textit{incorrect} problem that was proposed by Simmons in \cite{SM}: 

\textit{Given the initial value problem 
\begin{equation}\label{EulerEx}
y'=y^2+1,\;\;y(0)=0.
\end{equation}
 Use the Euler Method  to estimate numerically $y(2)$ with steps $h=$0.4,  $h=$0.2 and $h=$0.1 starting from the origin. Discuss your results.}

\begin{figure}[h]
\begin{center}
\includegraphics[width=0.5\textwidth]{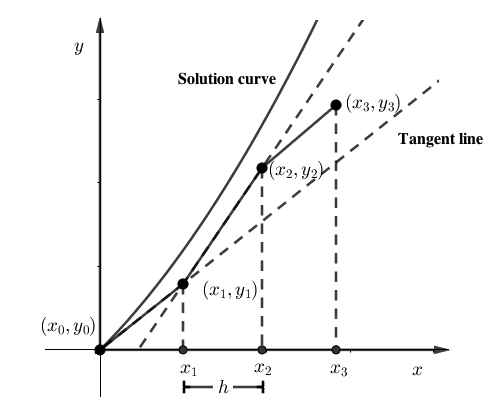}
\caption{The figure shows three iterations to find the value of $y_3$ that approximates the image $y(x_3)$ of the solution curve of a first order differential equation.}
\label{figure:Euler}  
\end{center}
\end{figure}

The problem given is an initial value problem of the form 
\begin{equation}\label{ODE}
\frac{dy}{dx}=f(x,y), \;\;y(x_0)=y_0.
\end{equation}
In order to find $y_k=y(x_k)$, the Euler method-- as also do other numerical methods to solve differential equations-- considers a sequence of points $\{x_k\}$ starting from the initial condition $x_0$ and separated by a step size $h$. It then uses the function $f(x,y)$ in Eqn. (\ref{ODE}) that represents the slope of the tangent line to the solution curve at any given point to find the next image in the form $y_{k+1}=y_k+f(x_k,y_k)h$ (see Fig. \ref{figure:Euler} ). Starting from $(x_0,y_0)$ several iterations then return the desired $y_k=y(x_k)$ value that approximates the image of the solution curve at $x_k$ \cite{Zill}.

Table \ref{EulerMeth} shows some results for $y(2)$  when performing some iterations of Euler's Method for the example from Eqn. (\ref{EulerEx}).

\begin{table}[h!]
\centering
\caption{The table shows in the last rows the approximate value of $y(2)$ obtained after 10 and 5 iterations of the Euler Method when using the step size  $h=0.2$ and $h=0.4$ respectively. }
\label{EulerMeth}
\begin{tabular}{ccc|ccc} \hline  \\[-.25cm]
\multicolumn{3}{c}{$h=0.2$}& \multicolumn{3}{c}{$h=0.4$}\\[.2cm]
$n$ & $x_n$ & $y_n$& $n$ & $x_n$ & $y_n$\\ [.2cm]
\hline   \\[-.2cm]
0& 0 & 0 & 0 &0  &  \\[.13cm]
1 & $0.2$ & $0.2$  & 1 & $0.4$ & $0.4$  \\[.13cm]
2& $0.4$ & $0.408$ & 2 & $0.8$  & $0.864$ \\[.13cm]
3 & $0.6$ & $0.6413$ & 3 & $1.2$ & $1.5626$ \\[.13cm]
4 & $0.8$ &  $0.9235$ & 4 & $1.6$ & $2.9393$ \\[.13cm]
5 & $1.0$ &$1.2941$ & 5 & $2$ & $6.795$ \\[.13cm]
6& $1.2$ & $1.8291$ &  &  &  \\[.13cm]
7 & $1.4$ & $2.6982$ &  &  &  \\[.13cm]
8 & $1.6$ & $4.3542$ &  &  &  \\[.13cm]
9 & $1.8$ & $8.3461$ &  &  &  \\[.13cm]
10 & $2$ &  $22.4778$ &  &  &  \\[.13cm]
\hline  
\end{tabular}
\end{table}

The values observed in Table \ref{EulerMeth} gives a starting point to think that something strange is happening in this example when using this method. We could ask ourselves: Why is there such a high variability in the approximate values of $y(2)$ when choosing different step sizes? The first approach to start understanding this, is to consider whether the initial value problem has a unique solution, since the variability could be due to the method finding different solutions, and hence $y(2)$ could be different for each solution. Analyzing the hypothesis of the Theorem of Existence and Uniqueness for first order initial value problems \cite{Zill}, one concludes easily that the initial value problem from Eqn. \ref{EulerEx} has a unique solution. 
Therefore the variation is likely due to some other cause.

A second reflection could be to look at the differential equation in Eqn. (\ref{EulerEx}) as an autonomous equation and examine the dynamics of its family of solutions. Since for the given equation $y'=y^2+1$ there are no real constant solutions, the solutions  could perfectly tend to infinity when approaching $x=2$. Hence students could want to justify the high variability of the solution $y$ at $x=2$ stating that the solution goes to infinity at $x=2$. 
  
 A third reflection could be to try to find an exact solution of the initial value problem, which when solving the differential equation for instance as a separable equation and using the initial condition gives us  $y=\tan(x).$ The student can note that the interval of definition of this solution is the  interval $\displaystyle ]-\frac{\pi}{2},\frac{\pi}{2}[$ and that the solution does not exist beyond $\pi/2\approx$ 1.57079. Thus, there is no point in asking for $y(2)$ because $x=2$ is not within the interval of definition of the solution. This reflection can lead to students realizing that the theorem of existence and uniqueness is a local theorem, and therefore only guarantees the existence of a solution in a neighborhood of the initial condition-- in this example in the neighborhood of $(0,0)$-- and that the interval of definition of the solution of a differential equation depends on the given initial condition.  
 
\begin{figure}[h!]
\begin{center}
\includegraphics[width=0.46\textwidth]{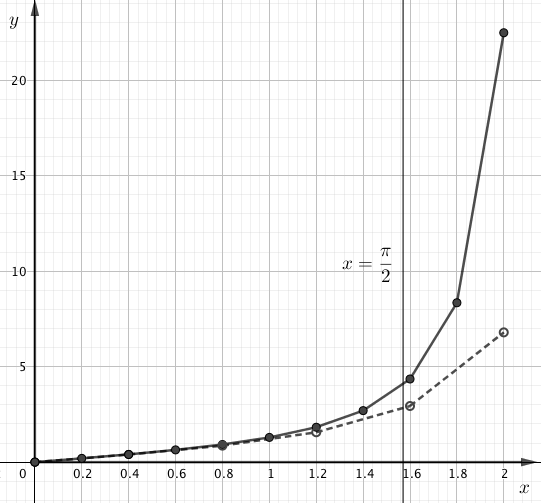}
\caption{The figure shows the approximate solution curve for the initial value problem given in Eqn. (\ref{EulerEx}). The solid line shows the curve using the Euler Method with step size $h=0.2$ and the dashed line using the step size $h=0.4$. }
\label{figure:2}
\end{center}
\end{figure}

 After the above reflections, a question that students should be confronted with and will need to answer is:  why when using the Euler method do we not realize that it does not make sense to ask for $y(2)$? If students are sufficiently persistent, they should reach to the conclusion this is the case since the value of each step size  $h$ is a rational number and hence starting from the origin, no matter how many steps of length $h$ we make, the point on the $x$-coordinate will never coincide with the irrational number $\pi/2$. Hence, the tangent line taken at the last partition point before $\pi/2$ will cross $\pi/2$ without noticing that it has already escaped the interval of definition of the solution as can be seen in Fig. \ref{figure:2} for the step sizes $h=0.2$ and $h=0.4$.

 
%

\subsection{A problem that makes us reflect on how temperature of an object decreases.} \label{Cooling}

In this subsection we present a problem that requires to model-- using differential equations-- how the temperature of an object changes with time, as well as to use calculus and algebra skills to analyze the model and find a solution. This problem leads to an interesting reflection regarding the differential equation model used that we will discuss further in Sec. \ref{Students}.  
Let us consider the following \textit{incorrect} problem:

 \textit{We have a pool with water at a constant temperature of $T_M.$ We place a thermometer into the pool, which initially shows $40^{\circ}C.$  After half a minute it marks $36^{\circ}C$ and after one minute $30^{\circ}C$.  What is the temperature of the water in the pool?}

Newton's Law of Cooling states that the temperature variation of an object in time $t$ is proportional to the difference between the temperature of the object and the temperature of the medium \cite{Zill}. In this case the temperature of the medium is the temperature of  the pool's water, which is assumed to have constant temperature $T_M$. In other words, if $T(t)$ represents the temperature of the thermometer at time $t,$ then 
\begin{equation}\label{ODECooling}
\frac{dT}{dt}=k(T-T_M).
\end{equation}
By solving the differential equation for instance as a separable equation we obtain 
\begin{equation}\label{TSol}
T(t)=T_M+(40-T_M)e^{kt}.
\end{equation}
Replacing the data points provided in the problem, we obtain the following system 

\begin{equation}
\left.\begin{array}{ccc} T_M+(40-T_M)e^{k/2}& = & 36 \\T_M+(40-T_M)e^{k} & = & 30\end{array}\right\}.
\end{equation}
If we solve for the exponential, we obtain 
\begin{equation}\label{expk}
e^{k/2}= \frac{36-T_M}{40-T_M}\;\;\mbox{and}\;\;e^{k}  =  \frac{30-T_M}{40-T_M}. 
\end{equation}
This leads to the following equation which must be satisfied by $T_M$:

\begin{equation}\label{eqn:T_M}
(36-T_M)^2=(30-T_M)(40-T_M).
\end{equation}
 After manipulation, the term $T_M^2$ is canceled, which results in $T_M=48^{\circ}C.$ In other words, by placing a thermometer showing an initial temperature of $40^{\circ}C$ in a pool at $48^{\circ}C$, the thermometer lowers the water's temperature, which is something that goes against common sense. Hence one realizes that something went wrong and the problem is an \textit{incorrect} problem. This motivates a discussion and a deeper understanding of the modeling aspects of the problem when we ask students: Where is the mistake in the reasoning and/or solution method of this exercise?

\section{Exposing students to an \textit{incorrect} problem in Differential Equations}\label{Students}

The problem that we presented in Subsec. \ref{Cooling} arose from an exam proposal for a basic course of Differential Equations for engineering students and was intended to be a relatively basic application of Newton's Law of  Cooling. After revising the proposal, the professors lecturing the course realized that this is an \textit{incorrect} problem and how many interesting questions it provokes that we could ask our students to think about.
Hence we decided to design a homework assignment to 168 undergraduate engineering students enrolled in the course Differential Equations. In the assignment, we firstly asked for the solution of a ``correct" problem and secondly asked to analyze and discuss an \textit{incorrect} problem that led to contradictions. Students worked on it for two weeks after which they had to turn in their assignment. Then, the assignments were graded and students' performance was analyzed. In the following Subsec. \ref{ProbCooling} we describe the assignment that was given to the students and in Subsec. \ref{Analysis} we discuss students' performance in both, the ``correct" and the \textit{incorrect} part of the problem that was assigned.  

\subsection{Students' assignment: the ``correct" problem and the \textit{incorrect} problem}\label{ProbCooling}

Students were confronted with two problems in their homework assignment. The first one is a ``correct" problem that can be solved directly since the problem is well posed and does not provide any contradictions in the solution process. Students in this course are generally used to these types of ``correct" problems. The problem is the following:\\

\textit{A thermometer showing $40^{\circ}$C is placed inside a pool filled with water that has a constant temperature of $T_M$. An observer notices that the thermometer shows $34^{\circ}$C after $1/2$ minute and $30^{\circ}$C after 1 minute.
\begin{enumerate}
 \item[a)] Determine the differential equation which models the dynamic of the temperature shown by the thermometer while it is in the pool.
 \item[b)] Using the equation from the previous item and the given data points: what is the temperature $T_M$ of the pool's water?\\
\end{enumerate}}
Using the model and solution method of Subsec. \ref{Cooling}, without further difficulties one obtains that $T_M=22^{\circ}C$, which makes perfect sense and solves the problem without further analysis. \\

The second problem that we assigned to students is the \textit{incorrect} problem, which  we also presented in Subsec. \ref{Cooling}, it states as follows: \\

\textit{A thermometer showing $40^{\circ}$C is placed inside a pool filled with water that has a constant temperature of $T_M$. An observer notices that the thermometer shows $36^{\circ}$C  after $1/2$ minute and $30^{\circ}$C after 1 minute.
Discuss and answer the following questions:
\begin{enumerate}
\item[a)] What is the temperature $T_M$, of the pool's water?
\item[b)] Think about how exactly the temperature of the thermometer within the pool changes with time. Discuss the contradiction that arises given the value of  $T_M$ obtained in a).
\item[c)] Plot the solution curve of the problem for the temperature $T_M$ obtained in a).
\item[d)] What happens to the value of $T_M$ when changing $T(1/2) = 36^{\circ}$C to $T(1/2) = 35^{\circ}$C ?
\item[e)] What can you say about the value of $k$ so that the differential equation models correctly a real situation such as the one observed in the pool?
\item[f)] For what range $T(1/2)$ values will the differential equation model the real situation correctly, given the presented data?
\end{enumerate}}

 As discussed in Subsec. \ref{Cooling}, we denote this problem \textit{incorrect} problem since when solving item $a)$ we obtain for the temperature of the water in the pool $T_M=48^{\circ}C$, which is a result that does not make sense in the given context, as we asked our students to discuss in item $b)$. 
 Next,  in Subsec. \ref{Analysis} we discuss the answers to items $b)- f)$ and how students performed in each of the items. 
\subsection{Answers to the \textit{incorrect} problem and students' performance}\label{Analysis}

We collected and analyzed the answers of the 168 undergraduate engineering students across 7 course sections that were enrolled in the Differential Equations course during one semester. Table \ref{Results} summarizes the results on student performance, showing the percentage students that provided a correct solution for each of the items of the assignment described in Subsec. \ref{ProbCooling}. 

\begin{table}[h!]
\centering
\caption{The table shows the percentage of correct answers of the students' assignment for the ``correct" problem and for the different items $a) -c)$, $d)$, $e)$ and $f)$ of the \textit{incorrect} problem. }
\label{Results}
\begin{tabular}{l|c|cccc} \hline  \\[-.25cm]
Type of& Correct & \multicolumn{4}{c}{Incorrect}\\[.2cm]
Problem&  & \\
&$a)$, $b)$ & $a)-c)$& $d)$ & $e)$ & $f)$\\ [.2cm]
\hline   \\[-.2cm]
 $\%$ of correct & 91$\%$ & 62$\%$ & 10$\%$ & 21$\%$ & 0$\%$\\
 answers &  & &  &  &  \\[.13cm]
\hline  
\end{tabular}
\end{table}
\FloatBarrier


The results of analyzing students' performance in solving the ``correct" problem show that most of the students were able to answer that problem. In fact,  91$\%$ of the students answered this question correctly. This shows that students did not have any difficulties in solving a problem that is well posed: they new the model that needed to be used and had the necessary skills to solve the differential equation and the algebraic skills to obtain the value for $T_M$ using the data points given in the problem.   



On the contrary, the results when analyzing the performance of the same students in solving the \textit{incorrect} problem show the lack of deep understanding students have in what really matters when modeling a cooling process:
Item $a)$ of the \textit{incorrect} problem could be solved directly using the typical solution method as in Subsec. \ref{Cooling}, obtaining $T_M=48$. As mentioned before, this result for $T_M$ does not make sense in the given context. We asked students in item $b)$ to discuss this contradiction. Observe that according to the data points given in the problem, a cooling process should be happening (the temperature shown by the thermometer decreases) and hence the temperature of the pool's water has to be at most $30^{\circ}$ to allow it. The mathematical skills to solve items $a)-c)$ were similar to the ones for the ``correct" problem, the difference now was that the result did not make sense. The number of students answering these items of the \textit{incorrect} problem correctly decreased to 62 $\%$ compared to the ``correct" problem that was answered by 91$\%$ of the students correctly  (see Table \ref{Results}).

\begin{figure}[h!]
\begin{center}
\includegraphics[width=0.51\textwidth]{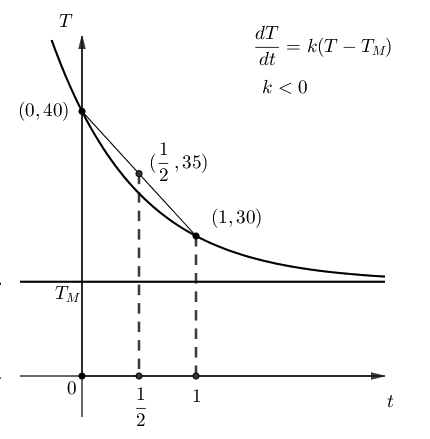}
\caption{The figure shows an exponential solution curve for the equation $\frac{dT}{dt}=k(T-T_M)$, $k<0$, for a model with data points $T(0)=40$, $T(1/2)<35$ and $T(1)=30$. It also depicts the secant line between the points $(0,40)$ and $(1,30)$ and illustrates that those two points are colineal with the point $(1/2,35)$.  }
\label{figure:7}
\end{center}
\end{figure}

When analyzing the answer for item $d)$ we observed that a small number of  students realized that since the dynamics of the temperature is modeled by a decreasing exponential (see Equation (\ref{TSol})), which is a convex function, and any image on the graph of this function between the points $(0,40)$ and $(1,30)$ will be smaller than the corresponding image on the graph of the line unifying these two points (see Fig. \ref{figure:7}). Hence, since the point $(1/2,35)$ is colinear with the points $(0,40)$ and $(1,30)$, the exponential curve can not be fitted to contain all three points and hence  the value of $T_M$ can not be obtained in particular with the solution method from  Subsec. \ref{Cooling}. Students had difficulty analyzing why the solution method failed, in fact, only 10 $\%$ of the students justified correctly why $T_M$ could not be found.

In item $e)$ we expected students to think about under which conditions the model actually describes the cooling process that we observe without obtaining a contradiction, and specifically asked for the value of the constant $k$. Observe that in a cooling process the temperature of the medium is less than the temperature of the object at any moment in time, and hence in Eqn. (\ref{ODECooling}) we have that $T-T_M$ is positive and hence the constant $k$ must be negative in order for the temperature to be decreasing ($dT/dt<0$) and for the model to describe well the cooling process.  This does not happen here, since if we solve for $k$ from Eqn. (\ref{expk}) when replacing by $T_M=48^{\circ}C$, we obtain $k=0.81$, and hence the solution method produces a contradiction with the model. Only a small group of 21 $\%$ of the students was able to explain the reason as to why the constant should be negative and that the solution method produces a contradiction because it produces a positive $k$ value.

\begin{figure}[h!]
\begin{center}
\includegraphics[width=0.47\textwidth]{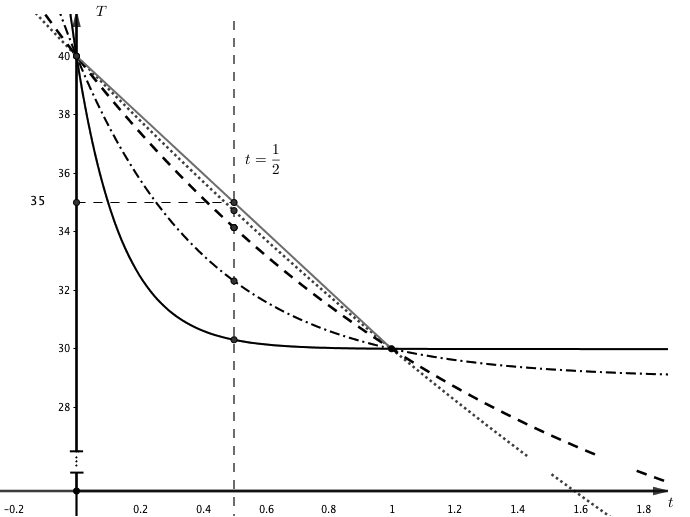}
\caption{solid for $ T_M=29.99$, $30.31$, point dashed  $T_M=29 $, $T(1/2)=32.32$, dashed for $T_M=20 $, $T(1/2)=34.14$, point for $T_M=-10 $ $T(1/2)=34.72$, light grey the secant line}
\label{figure:9}
\end{center}
\end{figure}

There were  no students who analyzed correctly item $f)$. Observe that if we replace the data points $T(1)=30$ and $T(1/2)=c$ in the solution process described in Subsec. \ref{Cooling} and then solve for $T_M$ in Eqn. (\ref{eqn:T_M}) we obtain $T_M$ as a function of $T(1/2)=c$ in the following way  \begin{equation}\label{TMc}
 T_M=\frac{1200-c^2}{70-2c}.
\end{equation}
Observe that $\displaystyle \lim_{c\rightarrow 35^-}T_M=-\infty$ and thence the real situation gets modeled correctly only for $c$ values for which $T_M$ attains temperature values larger then the absolute minimum ($-273.15^{\circ}C$). In other words, what cannot happen is that $T_M<-273.15$ and using the expression in Eqn. (\ref{TMc}) we obtain that this happens if  $T(1/2)=c$ is in the range  $34.9594<c<35$. So, for the model to describe a real situation, the data point $T(1/2)=c$ has to be in the range $30<c\leq 34.9594$. It can be seen from Fig, \ref{figure:9} what happens, the figure  depicts solution curves $T(t)=T_M+(40-T_M)e^{kt}$, $k<0$, for different data points $T(1/2)$. 
Observe that $\displaystyle \lim_{t\rightarrow \infty}T(t)=T_M$ and  observe from the figure that the closer $T(1/2)$ is to $35$, the lesser the limit $T_M$ is.  Hence we conclude that it is important to notice that the value of the data point $T(1/2)$ has a threshold below $35^{\circ}C$-- which here is $34.9594^{\circ}C$--  such that above this threshold the model does not make sense.  Questions like whether there are materials that do not behave exponentially as temperature decreases did not appear in the student's analyses and would be interesting to explore further.


 \section{Discussion}\label{Concl}

In this work we showed examples of \textit{incorrect} problems and displayed possible ways these could be used to motivate critical mathematical thinking in students. In particular, in Sec. \ref{ExInc} we described \textit{incorrect} problems that illustrated how important hypothesis are when modeling a problem or using a technique and how important it is to understand the restrictions of methods (numerical methods). We exposed a group of students to an \textit{incorrect} problem in Differential Equations which gave us important insight on how students react to a problem, which gives rise to a contradiction in their mathematical solution method or in their interpretation in a real life context. 

We observed several important aspects when confronting our Differential Equations students with an \textit{incorrect} problem: 
First, in order use \textit{incorrect} problems as a teaching opportunity, we think it is very important to first engage students with these types of \textit{incorrect} problems in a way they reach a level of natural curiosity regarding the questions posed by the problem. This is very important, since nobody makes discoveries without first being interested in the topic. Second, we think that \textit{incorrect} problems that are properly posed and guided can lead to a better analysis and understanding of the hypotheses and methods used to solve a problem. For many students reaching this analytic level is complicated due to the weaknesses and mathematical qualifications prior to the course, which need to be consolidated in order to tackle these topics. For instance, those who have not heard of convexity in Subsec. \ref{Analysis} are clearly at a disadvantage for properly performing a stronger analysis of the causes of the contradictions in the model.
Third, we must be self-critical with regards to having accustomed our students to respond rather direct ``correct" questions, which is why some fail when analyzing situations that escape the established boundaries. 

In history, people questioning their learning, rebelling against the established models are the ones that have been capable of great scientific revolutions, as shown by Einstein or Galileo, who questioned their academic realities and discovered new representations of our universe (see \cite{EinsteinPais,EinsteinIsaacson,EinsteinGoberoff,GalileoGeneralDiscoveries} for great literature on the life of individuals that changed the vision of the universe). Thus, we believe it to be necessary--  after the initial knowledge of a topic is learned-- to present students also with questions that require ``out of the box" thinking, additional analysis and even contradictions, allowing notable students in our course to discover new viewpoints and transforming themselves, viewing their knowledge from a different perspective. It is important to develop  reflections on mathematical topics that infuse life to our science and show it as novel, defying and attractive to all who we invite to enter the wonderful world of discovery and mathematics.

\section*{Acknowledgements}
We thank our colleagues  from the Mathematics Department of the Universidad Adolfo Ib\'a\~nes that taught Differential Equations that semester for giving us important feedback on this work. We also thank Clemens Vogt Geisse for his valuable contribution in translating this manuscript from its Spanish version.





\bibliographystyle{apa}

\end{document}